%
%

\magnification= \magstep1

\tolerance=300
\pretolerance=200
\hfuzz=1pt
\vfuzz=1pt
\parindent=35pt
\mathsurround=1pt
\parskip=1pt plus .25pt minus .25pt
\normallineskiplimit=.99pt

\hsize=5.8 true in
\vsize=9.5 true in

\hoffset=210 true mm             
\voffset=297 true mm             
\advance\hoffset by -\hsize      
\divide\hoffset  by 2            
\advance\hoffset by -1 true in   

\advance\voffset by -\vsize      
\divide\voffset  by 2            
\advance\voffset by -1.5 true in   

\font\eightrm=cmr8

\font\bfone=cmbx10 scaled\magstep1

\font\eightbf=cmbx8

\def\nin{\noindent}



\def\rightheadline{\hfil}
\def\leftheadline{\hfil}
\def\firstheadline{\hfill}

\countdef\fpageno=202  \fpageno=1
\countdef\onepd=204    \onepd=0

\def\seite #1 {\pageno #1
               \headline={\ifnum\pageno=#1 \firstheadline
               \else\ifodd\pageno\rightheadline\else\leftheadline\fi\fi}}

\def\ppp{??}
\def\pppp{??--??}
\def\pages{\ifodd\onepd\ppp\else\pppp\fi}

\def\commence{
%
%
\seite {\the\fpageno} }

\def\datum{\ifcase\month\or January\or February \or March\or April\or May
\or June\or July\or August\or September\or October\or November
\or December\fi\space\number\day, \number\year}

\def\date{{\rm Version of~}\datum}

\def\title{\centerline{Title ??}}
\def\author{Author ??}

\def\thanks#1{\footnote{$^*$}{\eightrm#1}}
\def\firstpage{\vbox{\vskip2.5truecm}
   \centerline{\bfone\title}\vskip.8truecm%
   \centerline{\bf\author}\vskip.8truecm\rm}

\def\abstract #1{{\narrower{\noindent\eightbf Abstract.\quad \eightrm #1 }
\bigskip}}

\def\sectionheadline #1{\bigbreak\vskip-\lastskip\indent\vskip1truecm
                       \centerline{\bf #1}\nobreak\medskip\nobreak}

\newtoks\literat
\def\[#1 #2\par{\literat={#2\unskip.}%
  \hbox{\vtop{\hsize=.1\hsize\nin [#1]\hfill}
  \vtop{\hsize=.9\hsize\nin\the\literat}}\par\vskip.3\baselineskip}
\def\references{
\sectionheadline{\bf References}
\frenchspacing \entries\par}

\countdef\revised=100
\revised=0

\def\recdate{}

\def\address{Author: {\tt$\backslash$def$\backslash$address$\{$??$\}$}}
\def\addresstwo{}

\def\lastpage{\references%
              \parindent=0pt\eightrm \baselineskip=8pt
              \vskip.8truecm
\line{\vtop{\hsize=.45\hsize{\eightrm\parindent=0pt\baselineskip=8pt
            \nin\address}}\qquad 
      \vtop{\hsize=.45\hsize\nin{\eightrm\parindent=0pt
            \baselineskip=8pt\nin\addresstwo}}\hfill}
\vskip.8truecm\recdate%
}

\mathchardef\emptyset="001F 

\def\qed{{\unskip\nobreak\hfil\penalty50\hskip .001pt \hbox{}\nobreak\hfil
          \vrule height 1.2ex width 1.1ex depth -.1ex
           \parfillskip=0pt\finalhyphendemerits=0\medbreak}\rm}

\def\Proposition #1. {\bigbreak\vskip-\parskip\noindent{\bf Proposition #1.}
    \quad\it}
\def\Theorem #1. {\bigbreak\vskip-\parskip\noindent{\bf  Theorem #1.}
    \quad\it}
\def\Corollary #1. {\bigbreak\vskip-\parskip\nin{\bf Corollary #1.}
    \quad\it}
\def\Lemma #1. {\bigbreak\vskip-\parskip\noindent{\bf  Lemma #1.}\quad\it}
\def\Definition #1. {\rm\bigbreak\vskip-\parskip\noindent{\bf Definition #1.}
    \quad}
\def\Remark #1. {\rm\bigbreak\vskip-\parskip\noindent{\bf Remark #1.}\quad}
\def\Exercise #1. {\rm\bigbreak\vskip-\parskip\noindent{\bf Exercise #1.}
    \quad}
\def\Example #1. {\rm\bigbreak\vskip-\parskip\noindent{\bf Example #1.}\quad}
\def\Examples #1. {\rm\bigbreak\vskip-\parskip\noindent{\bf Examples #1.}\quad}
\def\Proof#1.{\rm\par\ifdim\lastskip<\bigskipamount\removelastskip\fi
    \smallskip\noindent {\bf Proof.}\quad}

\def\AmS{$\cal A\kern-.1667em\lower.5ex\hbox{$\cal M$}\kern-.075em S$}

\def\ref#1{\expandafter\edef\csname#1\endcsname}
\ref {Introduction}{1}
\ref {Topology}{2}
\ref {ConvexCrit}{Theorem\penalty 10000\ 2.2}
\ref {E1}{$(2.3)$}
\ref {SjaThm}{Theorem\penalty 10000\ 2.3}
\ref {Open}{Proposition\penalty 10000\ 2.4}
\ref {Conn}{Lemma\penalty 10000\ 2.5}
\ref {ConvexCrit2}{Theorem\penalty 10000\ 2.6}
\ref {CP1}{Theorem\penalty 10000\ 2.7}
\ref {CP2}{Corollary\penalty 10000\ 2.8}
\ref {Counter}{3}
\ref {E3}{$(3.1)$}
\ref {ConvExamples}{4}
\ref {Proj}{Proposition\penalty 10000\ 4.1}
\ref {Examp}{Theorem\penalty 10000\ 4.2}
\ref {LocalConv}{5}
\ref {SKL}{Theorem\penalty 10000\ 5.1}
\ref {Property}{6}
\ref {E2}{$(6.1)$}
\ref {CP3}{Corollary\penalty 10000\ 6.1}
\ref {fa0}{Corollary\penalty 10000\ 6.2}
\ref {References}{7}
\ref {GS}{[1]}
\ref {HeHu}{[2]}
\ref {HNP}{[3]}
\ref {KL}{[4]}
\ref {Kir}{[5]}
\ref {Kn}{[6]}
\ref {Le}{[7]}
\ref {Sja}{[8]}

\def\cite#1{\csname#1\endcsname}
\newcount\Abschnitt\Abschnitt0
\def\beginsection#1. #2 \par{\advance\Abschnitt1\SATZ1\GNo0%
\sectionheadline{\bf\number\Abschnitt. \ignorespaces#2}\noindent}
\def\Definition #1. {\rm\bigbreak\vskip-\parskip\noindent{\bf Definition.}
    \quad}
\def\Remark #1. {\rm\bigbreak\vskip-\parskip\noindent{\bf Remark.}\quad}

\newcount\SATZ\SATZ1

\def\proclaim #1. #2\par{\bigbreak\vskip-\parskip\noindent{\bf#1}\quad
{\it#2}\Par\bigbreak}

\def\Aussage#1{%
\expandafter\def\csname#1\endcsname##1.{\proclaim {#1~\number
\Abschnitt.\number\SATZ.\global\advance\SATZ1}.}}
\Aussage{Theorem}
\Aussage{Proposition}
\Aussage{Corollary}
\Aussage{Lemma}
\def\Bindestrich{\penalty10000-\hskip0pt}
\let\_=\Bindestrich
\def\Links#1{\llap{$\scriptstyle#1$}}
\def\Rechts#1{\rlap{$\scriptstyle#1$}}
\def\Par{\par}
\def\into{\hookrightarrow}
\def\pfeil{\rightarrow}
\def\untenPf{\downarrow}
\def\pf#1{\buildrel#1\over\rightarrow}
\def\Pf#1{\buildrel#1\over\longrightarrow}
\def\|#1|{\mathop{\rm#1}\nolimits}
\def\<{\langle}
\def\>{\rangle}
\def\P{{\bf P}}
\catcode`\@=11
\def\hex#1{\ifcase#1 0\or1\or2\or3\or4\or5\or6\or7\or8\or9\or A\or B\or
C\or D\or E\or F\else\message{Warnung: Setze hex#1=0}0\fi}
\def\fontdef#1:#2,#3,#4.{%
\alloc@8\fam\chardef\sixt@@n\FAM
\ifx!#2!\else\expandafter\font\csname text#1\endcsname=#2
\textfont\the\FAM=\csname text#1\endcsname\fi
\ifx!#3!\else\expandafter\font\csname script#1\endcsname=#3
\scriptfont\the\FAM=\csname script#1\endcsname\fi
\ifx!#4!\else\expandafter\font\csname scriptscript#1\endcsname=#4
\scriptscriptfont\the\FAM=\csname scriptscript#1\endcsname\fi
\expandafter\edef\csname #1\endcsname{\fam\the\FAM\csname text#1\endcsname}
\expandafter\edef\csname hex#1fam\endcsname{\hex\FAM}}
\catcode`\@=12 
\fontdef Fr:eufm10,eufm7,eufm5.
\def\fa{{\Fr a}}
\def\fh{{\Fr h}}
\def\fk{{\Fr k}}
\def\fl{{\Fr l}}
\def\ft{{\Fr t}}
\fontdef bbb:msbm10,msbm7,msbm5.
\fontdef msa:msam10,msam7,msam5.
\def\CC{{\bbb C}}
\def\RR{{\bbb R}}
\def\ZZ{{\bbb Z}}
\mathchardef\auf=\string"3\hexmsafam10
\def\Dq{{\overline{D}}}
\def\Vq{{\overline{V}}}
\def\Xq{{\overline{X}}}

\def\NewDisplay\eqno#1
#2$${#2\neweqno{#1}$$}
\everydisplay{\NewDisplay}
\newcount\GNo\GNo=0
\def\neweqno#1{
\global\advance\GNo1
\eqno{(\number\Abschnitt.\number\GNo)}}
\newcount\refzaehler\refzaehler0

\newif\ifxxx
\xxxfalse
\catcode`/=\active
\def/{\string/\ifxxx\hskip0pt\fi}
\def\TText#1{{\xxxtrue\tt#1}}
\def\L|Abk:#1|Sig:#2|Au:#3|Tit:#4|Zs:#5|Bd:#6|S:#7|J:#8|xxx:#9||{%
\global\advance\refzaehler1
\[\the\refzaehler{} \ifx-#3\relax\else{#3}: \fi
\ifx-#4\relax\else{\it #4\/}{\sfcode`.=3000.} \fi
\ifx-#5\relax\else{#5} \fi
\ifx-#6\relax\else{\bf #6} \fi
\ifx-#8\relax\else({#8})\fi
\ifx-#7\relax\else, {#7}\fi\ifx-#9\relax\else, \TText{#9}\fi}
\def\B|Abk:#1|Sig:#2|Au:#3|Tit:#4|Reihe:#5|Verlag:#6|Ort:#7|J:#8|xxx:#9||{%
\global\advance\refzaehler1
\[\the\refzaehler{} \ifx-#3\relax\else{#3}: \fi
\ifx-#4\relax\else{``#4.''} \fi
\ifx-#5\relax\else{(#5)} \fi
\ifx-#7\relax\else{#7:} \fi
\ifx-#6\relax\else{#6}\fi
\ifx-#8\relax\else{ #8}\fi
\ifx-#9\relax\else, \TText{#9}\fi}

\def\r#1{{\it\romannumeral#1)}}
\def\rs#1{{\it$\it\romannumeral#1'$)}}
\def\dsseq{\hbox{$\cap$\vrule}}

\def\date{{\rm Version of~}Dec. 13, 2001} 

\def\title{\centerline{Convexity of Hamiltonian manifolds}}

\def\author{Friedrich Knop\thanks{This work was supported by grants of the
NSA and the NSF.}}

\def\address{Department of Mathematics\par
Rutgers University\par
New Brunswick NJ 08903\par
USA\par
knop@math.rutgers.edu
}
\def\addresstwo{} 

\def\entries{%

\B|Abk:GS|Sig:GS|Au:Guillemin, V.; Sternberg, S.|Tit:Symplectic
techniques in physics|Reihe:-|Verlag:Cambridge University
Press|Ort:Cambridge|J:1984|xxx:-||

\L|Abk:HeHu|Sig:HH|Au:Heinzner, P.; Huckleberry, A.|Tit:K{\"a}hlerian
potentials and convexity properties of the moment map|%
Zs:Invent. Math.|Bd:126|S:65--84|J:1996|xxx:-||

\L|Abk:HNP|Sig:HNP|Au:Hilgert, J.; Neeb, K.-H.; Plank, W.|%
Tit:Symplectic convexity theorems and coadjoint orbits|%
Zs:Compos. Math.|Bd:94|S:129--180|J:1994|xxx:-||

\L|Abk:KL|Sig:KL|Au:Karshon, Y.; Lerman, E.|Tit:The centralizer of
invariant functions and division properties of the moment map|%
Zs:Illinois J. Math.|Bd:41|S:462--487|J:1997|xxx:dg-ga/9506008||

\L|Abk:Kir|Sig:Ki|Au:Kirwan, F.|Tit:Convexity properties of the moment
mapping III|Zs:Invent. math.|Bd:77|S:547-552|J:1984|xxx:-||

\L|Abk:Kn|Sig:Kn|Au:Knop, F.|Tit:Weyl groups of Hamiltonian manifolds,
I|Zs:Preprint|Bd:-|S:34 pages|J:1997|xxx:dg-ga/9712010||

\L|Abk:Le|Sig:Le|Au:Lerman, E.|Tit:Symplectic cuts|Zs:Math. Research
Lett.|Bd:2|S:247--258|J:1995|xxx:-||

\L|Abk:Sja|Sig:Sj|Au:Sjamaar, R.|Tit:Convexity properties of the moment map
re-examined|Zs:Adv. Math.|Bd:138|S:46--91|J:1998|xxx:dg-ga/9408001||

}

\commence
\firstpage

\abstract{We study point set topological properties of the moment
map. In particular, we introduce the notion of a convex Hamiltonian
manifold. This notion combines convexity of the momentum image and
connectedness of moment map fibers with a certain openness requirement
for the moment map. We show that convexity rules out many pathologies
for moment maps. Then we show that the most important classes of
Hamiltonian manifolds (e.g., unitary vector spaces, compact manifolds,
or cotangent bundles) are in fact convex. Moreover, we prove that
every Hamiltonian manifold is locally convex.}

\let\oldbracket=\[
\def\[#1,#2]{\overline{#1#2}}

\beginsection Introduction. Introduction

Let $K$ be a connected compact Lie group with Lie algebra $\fk$ and let
$M$ be a Hamiltonian $K$\_manifold, i.e., a symplectic $K$\_manifold
equipped with a moment map $\mu:M\pfeil\fk^*$. The purpose of this
note is to study certain point\_set topological properties of $\mu$.

Let $\ft\subseteq\fk$ be a Cartan subalgebra and $\ft^+\subseteq\ft^*$
a Weyl chamber. Since every $K$\_orbit of $\fk^*$ meets $\ft^+$ in
exactly one point, the image of $\mu$ is determined by
$\psi(M):=\mu(M)\cap\ft^+$. A celebrated theorem of Kirwan,
\cite{Kir}, states that if $M$ is compact and connected then

\item{\r1} {\it $\psi(M)$ is a convex set,} and
\item{\r2} {\it all fibers of $\mu$ are connected.}

\noindent On the other hand, many other Hamiltonian manifolds have
these two properties. A very important example of a non\_compact
Hamiltonian manifold is a finite\_dimensional unitary representation
of $K$. Here also \r1, \r2 hold.  Therefore, it seems worthwhile to
introduce a general concept which encompasses both the compact and the
unitary case.

In this paper we introduce a certain completeness property called {\it
convexity}. To define it, let $\psi:M\pfeil\ft^+$ be the map which
assigns to $x\in M$ the point of intersection of $K\mu(x)$ with
$\ft^+$. This is a $K$\_invariant continuous map. For any two points
$u,v\in\ft^+$ let $\[u,v]$ be the line segment joining $u$ with
$v$. Then we call $M$ {\it convex\/} if $\psi^{-1}(\[u,v])$ is
connected for all $u,v\in\psi(M)$.

It is pretty immediate from the definition that convexity of $M$
implies \r1 and \r2. Less obvious is the fact that convexity entrains
also a more subtle property, namely that

\item{\r3} {\it the map $\psi:M\pfeil\psi(M)$ is open.}

In section~\cite{Topology} we show that, conversely, \r1,
\r2, and \r3 are equivalent to convexity. To this end, we use
results of Sjamaar, \cite{Sja}, on the local structure of the moment
map. It is an open problem whether \r3 is really needed to imply
convexity. In any case, convexity is a useful concept since it
rules out much of the pathological behavior a moment map may have.
To illustrate this we give various examples of bad moment maps in
section~\cite{Counter}.

On the other hand, we prove that many Hamiltonian manifolds occurring
in applications are in fact convex. For example, we show that a
connected Hamiltonian manifold is convex whenever it is compact, a
complex algebraic K\"ahler variety, or a cotangent bundle.

Then, in section~\cite{LocalConv}, we show that every Hamiltonian
manifold is locally convex. This makes the concept of convexity
useful even for arbitrary Hamiltonian manifolds. This fact
conceptualizes some topological considerations in the paper \cite{KL}
by Karshon\_Lerman.

The convexity of $M$ has consequences far beyond the topological
properties \r1--\r3. A few are mentioned in this paper like the fact
that for convex $M$, the momentum image $\psi(M)$ is locally a convex
polytope. For a much deeper application see the paper \cite{Kn} where
we were able to identify all collective functions on $M$ provided $M$
is convex.

\medskip
\Remark. This paper is an extended version of the second section of the
preprint~\cite{Kn}.
\medskip
\noindent{\bf Acknowledgment.}\quad I would like to thank
Y.~Karshon, E.~Lerman, S.~Tolman, and Ch.~Woodward for useful
discussions on topics related to this paper. I also thank the
referee for pointing out a gap in the proof of \cite{Examp}\r4 in a
previous version of this paper.

\beginsection Topology. Convex Hamiltonian manifolds: the definition and a
criterion

Let $K$ be a connected, compact Lie group with Lie algebra $\fk$. A
{\it Hamiltonian $K$\_manifold\/} is a $K$\_manifold $M$ with a
$K$\_invariant symplectic form $\omega$ and with a {\it moment map},
i.e., a $K$\_equivariant map $\mu:M\pfeil\fk^*$ such that
$$\eqno{}
\langle\xi,d\mu(\eta)\rangle=\omega(\xi x,\eta)\quad
\hbox{for all }\xi\in\fk, x\in M, \eta\in T_x(M).
$$
Let $\ft\subseteq\fk$ be a Cartan subalgebra corresponding to a
maximal torus $T\subseteq K$. Since
$\ft$ has a unique $T$\_stable complement in $\fk$, we can regard
$\ft^*$ as a subspace of $\fk^*$. Let $\ft^+\subseteq\ft^*$ be a
Weyl chamber. The composition $\ft^+\into\fk^*\auf\fk^*/K$
is a homeomorphism. We use it to construct a
continuous map $\psi$ which makes the following diagram commutative:
$$\eqno{}
\matrix{
M&\pf\mu&\fk^*\cr
\Links\psi\untenPf&&\untenPf\cr
\ft^+&\pfeil&\fk^*/K\cr}
$$
In other words, $\psi$ is the unique map with
$K\mu(x)\cap\ft^+=\{\psi(x)\}$ for all $x\in M$.

For any two (not necessarily distinct) points $u,v\in\ft^*$ let
$\[u,v]$ be the line segment joining them. Observe, that a subset $C$
of $\ft^*$ is convex if $\[u,v]\cap C$ is connected for all $u,v\in
C$. This is just a slight reformulation of the classical definition
and motivates:

\Definition. A Hamiltonian $K$\_manifold is called {\it convex\/} if
$\psi^{-1}(\[u,v])$ is connected for all $u,v\in\psi(M)$.
\medskip\noindent
A first reformulation of the concept is:

\Proposition. A Hamiltonian $K$\_manifold $M$ is convex if and only if
$\psi^{-1}(B)$ is connected for every convex subset $B$ of $\ft^+$.

\Proof. Suppose there is a convex subset $B\subseteq\ft^+$ such that
$\psi^{-1}(B)$ is disconnected. Then $\psi^{-1}(B)$ is the disjoint
union of two non\_empty open subsets, say $U$ and $V$. Let
$u\in\psi(U)$ and $v\in\psi(V)$. Since $B$ is convex, we have
$\[u,v]\subseteq B$, hence
$X:=\psi^{-1}(\[u,v])\subseteq\psi^{-1}(B)$. Thus, $X=(X\cap
U)\cup(X\cap V)$ is the disjoint union of two non\_empty open subsets
which means that $M$ is not convex. This shows one direction, the
other is trivial.\qed

\noindent Putting $u=v$ in the definition, we see that for a convex
manifold $M$ the fibers of $\psi$ are connected. Moreover, for every
$u,v\in\psi(M)$ the intersection $\[u,v]\cap\psi(M)$ has to be
connected, i.e., $\psi(M)$ is convex. These two properties are not
enough to imply convexity, though. More precisely, we have the
following characterization:

\Theorem ConvexCrit. A Hamiltonian $K$\_manifold is convex if and only
if the following conditions are satisfied:
\item{\r1} The image $\psi(M)$ is convex.
\item{\r2} The fibers of $\psi$ are connected.
\item{\r3} The map $\psi:M\pfeil\psi(M)$ is open.\Par

\noindent{\it Proof of ``$\Leftarrow$'':} Assume that \r1 through
\r3 hold. If $X:=\psi^{-1}(\[u,v])$ is disconnected then there are
open subsets $U_1,U_2$ of $M$ such that $X$ is the disjoint union of
the non\_empty subsets $X_i:=X\cap U_i$.  Suppose
$w\in\psi(X_1)\cap\psi(X_2)$. Then $F:=\psi^{-1}(w)$ is the union of
the non\_empty disjoint sets $F\cap U_1$ and $F\cap U_2$. Since this
contradicts \r2, we have $\psi(X_1)\cap\psi(X_2)=\emptyset$. Furthermore
$\psi(X_1)\cup\psi(X_2)=\psi(X)=\[u,v]\cap\psi(M)=\[u,v]$, by
\r1. Finally,
\r3 implies that $\psi(X_i)=\[u,v]\cap\psi(U_i)$ is open in
$\[u,v]$ which contradicts the connectedness of $\[u,v]$. The proof of
the reverse direction is deferred to after \cite{Conn}.\qed

All properties which we considered so far make sense for any
continuous map $\psi$ from a topological space $M$ to a convex subset
of a real vector space and we have shown that the following
implications hold:
$$\eqno{1}
\r1\wedge\r2\wedge\r3\quad\Rightarrow\quad\psi\hbox{ convex}
\quad\Rightarrow\quad\r1\wedge\r2
$$
Neither of these implications is reversible in this generality: For
the first arrow, let $\psi:M\pfeil\RR^2$ be the (real) blow up of
$\RR^2$ in the origin, i.e.,
$M=\{((x_0,x_1),[y_0:y_1])\in\RR^2\times\P^1(\RR)\mid
x_0y_1=x_1y_0\}$. Consider the inverse image $X$ of a line segment
$\[u,v]\subseteq\RR^2$. If $\[u,v]$ does not contain the origin then
$X\pf\sim\[u,v]$. Otherwise, $X$ is the union of the exceptional fiber
$\psi^{-1}(0,0)\cong S^1$ and an interval meeting in one point. In
either case, $X$ is connected, hence $\psi$ is convex. On the other
hand, the image of any sufficiently small neighborhood of a point in
the exceptional fiber is contained in a small cone and contains the
origin. Such a set can't be open, hence \r3 does not hold. For the
second arrow let $M$ be the same as above but with one point in the
exceptional fiber removed. Then $\psi$ is still surjective and all
fibers are connected. But there is a line segment $\[u,v]$ whose
preimage is disconnected.

This shows that we have to use some special properties of moment
maps. The dual of the Cartan algebra $\ft^*$ contains a canonical
lattice $\Gamma$ namely the differentials of all homomorphisms
$T\pfeil\RR/\ZZ$. A {\it rational homogeneous cone\/} is a subset
of $\ft^*$ which is of the form $\sum_{i=1}^N\RR_{\ge0}\gamma_i$ where
$\gamma_1,\ldots,\gamma_N\in\Gamma$. A {\it rational cone\/} is a
translate $u+C$ of a homogeneous rational cone $C$ by a vector $u\in
\ft^*$. In this case we say that $u$ is a vertex of $u+C$. Note that
$u+C$ can have many vertices namely all points in $u+(C\cap-C)$. The
following theorem of Sjamaar contains all the special properties of
the moment map which we are going to need.

\Theorem SjaThm. {\rm(\cite{Sja}~Thm.~6.5)} Let $M$ be a Hamiltonian
$K$\_manifold. Then for every orbit $Kx\subseteq M$ there is a unique
rational cone $C_x\subseteq\ft^*$ with vertex $\psi(x)$ such that:
\item{\r1} There exist an arbitrarily small $K$\_stable neighborhood $U$
of $Kx$ such that $\psi(U)$ is a neighborhood of $\psi(x)$ in $C_x$.
\item{\r2} For $u\in\ft^+$ let $x$ and $y$ be in the same connected
component of $\psi^{-1}(u)$. Then $C_x=C_y$.\Par

\noindent
A first application of this theorem is a general openness property of $\psi$:

\Proposition Open. Let $M$ be a Hamiltonian $K$\_manifold such that all
fibers of $\psi:M\pfeil\ft^+$ are connected. Let $U\subseteq M$ be
open. Then also $\psi^{-1}\psi(U)$ is open.

\Proof. Let $U'$ be the union of all translates $kU$, $k\in K$. Then
$U'$ is open with $\psi(U')=\psi(U)$. Thus, we may replace $U$ by $U'$
and assume that $U$ is $K$\_stable. Let $y\in\psi^{-1}\psi(U)$ and
$u:=\psi(y)$. Then there is $x\in U$ such that also $\psi(x)=u$. Since
$\psi^{-1}(u)$ is connected, \cite{SjaThm}\r2 implies that
$C_x=C_y=:C$. By part \r1 of that theorem there are open
neighborhoods $U_x$, $U_y$ of $x$, $y$, respectively, such that
$\psi(U_x)$ and $\psi(U_y)$ are neighborhoods of $u$ in $C$. Hence
$(\psi|_{U_y})^{-1}(\psi(U_x))$ is a neighborhood of $y$ which is
contained in $\psi^{-1}\psi(U)$. This proves the assertion.\qed

\noindent We can now derive an openness criterion for $\psi$:

\Lemma Conn. Let $M$ be a Hamiltonian $K$\_manifold such that all
fibers of $\psi:M\pfeil\ft^+$ are connected. Assume moreover that
every $u\in\psi(M)$ has an arbitrarily small neighborhood $B$ such that
$\psi^{-1}(B)$ is connected. Then $\psi:M\pfeil\psi(M)$ is an open
map.

\Proof. Let $x\in M$ and $U$ be an open neighborhood of $x$. We have to
show that $\psi(U)$ is a neighborhood of $u:=\psi(x)$ in
$\psi(M)$. For this we may assume that $\psi(U)$ is a neighborhood of
$u$ in $C_x$. Let $B$ be a neighborhood of $u$ in $\ft^+$ such that
$B\cap C_x\subseteq\psi(U)$ and such that $\psi^{-1}(B)$ is
connected. Then $V_1:=\psi^{-1}\psi(U)$ is open in $M$ by
\cite{Open}. Clearly, also $V_2:=\psi^{-1}(\ft^*\setminus C_x)$ is
open in $M$. Moreover, $V_1$ and $V_2$ are disjoint and cover
$\psi^{-1}(B)$. Connectivity implies $\psi^{-1}(B)\subseteq V_1$,
i.e., $B\cap\psi(M)\subseteq\psi(U)$ which proves the assertion.\qed

\noindent Now we can complete the proof of \cite{ConvexCrit}:
\medskip
\noindent{\it Proof of ``$\Rightarrow$'':} Assume that $M$ is
convex. Then \r1 and \r2 clearly hold. Let $u\in\psi(M)$ and
$B\subseteq\ft^+$ a convex neighborhood of $x$. Since also $\psi(M)$
is convex we have $\[u,v]\subseteq B\cap\psi(M)$ for every $v\in
B\cap\psi(M)$. By assumption, $\psi^{-1}(\[u,v])$ is connected. This
implies that $\psi^{-1}(B)$ is connected. Thus \r3 holds by
\cite{Conn}.\qed

One can express properties \r1--\r3 of \cite{ConvexCrit}
purely in terms of the moment map $\mu$:

\Theorem ConvexCrit2. A Hamiltonian $K$\_manifold is convex if and only
if the following conditions are satisfied:
\item{\rs1} The intersection $\mu(M)\cap\ft^+$ is convex.
\item{\rs2} The fibers of $\mu$ are connected.
\item{\rs3} Whenever $U\subseteq M$ is a {\rm $K$\_invariant\/} open
subset then $\mu(U)$ is open in $\mu(M)$.\Par

\Proof. By construction we have $\psi(M)=\mu(M)\cap\ft^+$, whence
\r1$\Leftrightarrow$\rs1. For $x\in M$ let $y:=\mu(x)$ and
$u:=\psi(x)$. Then there is a $K$\_equivariant map $\psi^{-1}(u)\pfeil
Ky$, hence $\psi^{-1}(u)=K\times^L\mu^{-1}(y)$ where
$L:=K_{\mu(x)}$. The equivalence of \r2 and \rs2 now follows
{}from the fact that $K$ and $L$ are connected. Finally, let $U\subseteq
M$ be open and $U'$ the union of all translates $kU$. Since
$\psi(U')=\psi(U)$ it suffices for checking \r3 that $\psi(U)$ is
open in $\psi(M)$ for all $K$\_invariant $U$. But then
\r3$\Leftrightarrow$\rs3 follows from the fact that
$\fk^*\pfeil\ft^+$ is the quotient map by $K$.\qed

\Remark. In general, one cannot expect $\mu$ to be an open map. Take
for example $M=T^*(K/H)$ where $K=SU(2)$ and $H=U(1)$, the maximal
torus. Then $M=K\times^H\fh^\perp$ where $\fh^\perp\subseteq\fk^*$ is
the annihilator of the Lie algebra of $H$. Then
$\mu(M)=\|ad|K\cdot\fh^\perp=\fk^*$. On the other hand, let
$V_0\subset K/H$ be a small open subset, $V$ its preimage in $K$ and
$U:=V\times^H\ft^\perp$. Then $U$ is open in $M$ but its image
$\mu(U)=V\fh^\perp$ is a conical {\it proper\/} subset of $\fk^*$
containing $0$. Thus, it is not open. 
\medskip

We conclude this section with two general properties of convex Hamiltonian
manifolds.  First we show that the cone $C_x$ can be recovered from
$\psi(M)$.

\Theorem CP1. Let $M$ be a convex Hamiltonian $K$\_manifold and $x\in
M$. Then $C_x$ is the smallest cone containing $\psi(M)$ and having
$\psi(x)$ as a vertex. Moreover, $\psi(M)$ forms a neighborhood of
$\psi(x)$ in $C_x$.

\Proof. Denote this smallest cone by $C$. By \cite{SjaThm}, there is a
neighborhood $U$ of $x$ such that $\psi(U)$ is a neighborhood of
$u:=\psi(x)$ in $C_x$. This shows that $C_x\subseteq C$. For the converse
we just have to show that $\psi(M)$ is contained in $C_x$.

Let $v\in\psi(M)$. Since $\psi(U)$ is open in $\psi(M)$ and since
$\[u,v]\subseteq\psi(M)$ also $\[u,v]\cap\psi(U)$ is open in
$\[u,v]$. This implies $\[u,v]\subseteq C_x$, thus $v\in C_x$.\qed

\Corollary CP2. Let $M$ be a convex Hamiltonian $K$\_manifold. Then
the image $\psi(M)$ is locally a polyhedral cone and, in particular,
locally closed and semi\_analytic in $\ft^*$.

\Remark. See Example~2 of the next section for a connected Hamiltonian
manifold with $\psi(M)\not\subseteq C_x$ and Example~3 where $\psi(M)$
is not locally closed.

\beginsection Counter. Pathological examples of moment maps

One can reformulate \cite{Open} as follows: If $\psi$ has
connected fibers then one can factorize it as
$$\eqno{3}
M\Pf\alpha S\Pf\beta\ft^+
$$
where $\alpha$ is open, surjective and $\beta$ is injective. In
fact, set-theoretically we let $S=\psi(M)$, $\alpha=\psi$, and
$\beta$ to be the natural inclusion of $S$ in $\ft^+$. We equip $S$ with the
topology that $V\subseteq S$ is open if and only if $\psi^{-1}(V)$ is
open in $M$.

How important is it that $\psi$ has connected fibers? There is a
straightforward generalization of the factorization \cite{E3}: let
$S:=M/\sim$ where $x\sim y$ if $x$ and $y$ are in the same connected
component of some fiber of $\psi$. We give $S$ the quotient
topology. Then the proof of \cite{Open} shows that $\alpha$ is an open
map. Moreover, the map $\beta$ is continuous with discrete
fibers.

\Example 1. Consider the above construction in the
following setting. Let $K=T$ be a torus and $M_0$ any Hamiltonian
$T$\_manifold with moment map $\mu_0$. Let $V$ be a manifold and
$\gamma:V\pfeil\ft^*$ a smooth map which is everywhere a local
isomorphism. Now consider the fiber product
$M=M_0\times_{\ft^*}V$. Since it is locally isomorphic to $M_0$, it is
also Hamiltonian. Its moment map is simply the projection to $V$
composed with $\gamma$. Thus, if $\mu_0$ has connected fibers then
$S=M/\sim=\gamma^{-1}\mu_0(M_0)$. This way, one can construct
Hamiltonian manifolds with arbitrarily disconnected fibers. A concrete
example is the following: take $T=(S^1)^2$,
$M_0:=T^*(T)=T\times\ft^*$, and identify $\ft^*\cong\RR^2$ with
$\CC$. Let $S=\CC\setminus\{0,{2\over3}\}$ and
$\gamma:S\pfeil\CC=\ft^*:z\mapsto z^3-z^2$. Then $\gamma$ is a
surjective unramified map. Since $\gamma$ is an open map, we have
constructed a Hamiltonian manifold $M=T\times S$ such that $\psi$ is
open with convex image but such that most fibers are disconnected.

\Example 2. Using the same technique as in
Example 1, one can also construct a Hamiltonian manifold such that
$\beta:S\pfeil\ft^*$ is injective but not open. For this take again
$K=(S^1)^2$, $M_0=T\times\ft^*$ and identify $\ft^*$ with $\CC$. Let
for some $0<\epsilon<\pi$ let $S'$ be the set of all non\_zero complex
numbers with $-\epsilon<\|Arg|z<\pi$ (i.e., $S'$ is a bit larger than
the upper half\_plane). Let $\gamma':S'\pfeil\CC=\ft^*:z\mapsto z^2$
and $M'=T\times S'$. Now we apply Lerman's symplectic cut technique,
\cite{Le}, to the preimage $Y$ of the positive real half line which is
isomorphic to $S^1\times S^1\times\RR$: one can cut away the preimage
of the sector $-\epsilon<\|Arg|z<0$ and replace $Y$ by $Y/(1\times
S^1)\cong S^1\times\RR$. Then one gets a new Hamiltonian manifold $M$
such that $S$ is the set of all $z\in\CC$, $z\ne0$ with
$0\le\|Arg|z<\pi$, i.e., $S$ is the upper half\_plane together with
the positive real half\_line. The map $\gamma:S\pfeil\CC:z\mapsto z^2$
is injective with image $\CC\setminus\{0\}$ but not open: no
neighborhood of the positive half\_line maps to a neighborhood of the
image point.  This is basically the Example~3.10 of \cite{KL}
and which is avoided by our notion of convexity.

\Example 3. This example is the same as the
preceding one but we remove the preimage of
$[1,\infty)\times(0,\infty)\subseteq\ft^*$. Then $\psi(M)$ is not
locally closed near the point $(1,0)$.

\Example 4. Unfortunately, $S$ will
not in general be Hausdorff. To show this, let $K=S^1$ and $M_0:=T^*(S^1)\times
V=S^1\times\RR\times V$ where $V$ is a non\_zero symplectic vector
space. The fibers of the moment map are $S^1\times V$. Now let
$H\subseteq V$ be a hyperplane and let $M$ be $M_0$ where we removed
$S^1\times H$ from one of the fibers of the moment map. Then one of
the fibers of the moment map gets disconnected and $S=M/\sim$ is a
real line with one of its points doubled. In particular, it is not
Hausdorff.

\Example 5. In the preceding example one could
remedy the situation by replacing ``$\sim$'' by a coarser equivalence
relation. More precisely, put $S:=M/\approx$ where $x\approx y$ if
$\psi(x)=\psi(y)$ and $C_x=C_y$. Then one can show that $\alpha$ is
still open and that $\beta$ has discrete fibers. But also in this
definition, $S$ might not be Hausdorff. To construct an example
consider the two Hamiltonian $S^1$\_manifolds
$M_1:=T^*(S_1)\times\RR^2=S^1\times\RR\times\RR^2$, and
$M_2:=\CC\times\RR^2$. Here $K=S^1$ acts in both cases on the first
factor and $\RR^2$ has the standard symplectic structure $dx\wedge
dy$. The moment maps are $\mu(\alpha,t,x,y)=t$ and
$\mu(z,x',y')=|z|^2$, respectively.  Consider the open subsets
$M_1^0:=\mu_1^{-1}(\RR^{>0})=\{t>0\}$ and
$M_2^0:=\mu_2^{-1}(\RR_{>0})=\{z\ne0\}$. They are isomorphic as
Hamiltonian manifolds, the isomorphism being $z=\sqrt{t}\|exp|(2\pi
i\alpha)$. Now we twist this isomorphism by the automorphism of
$M_1^0$
$$\eqno{}
(\alpha,t,x,y)\mapsto (\alpha-t^{-2}y,t,x+t^{-1},y)
$$
and glue $M_1$ and $M_2$ along $M_1^0\pf\sim M_2^0$. It is easily seen
that the graph of the gluing isomorphism is closed in $M_1\times M_2$
(on the graph holds $x'=x+t^{-1}$ and $|z|^2=t$. Thus, on its closure
we get $t(x'-x)=1$ and $|z|^2=t$ which implies $t>0$ and $z\ne0$,
i.e., the graph is closed). Therefore, the resulting manifold $M$ is
Hausdorff. This way we constructed a Hamiltonian $S^1$\_manifold for
which all fibers of the moment map but one are connected ($\cong
S^1\times\RR^2$) and the zero fiber is disconnected ($\cong
(S^1\times\RR^2)\mathop{\dot\cup}\RR^2$). Moreover, for $x$ in the
first component we have $C_x=\RR$ while in the second holds
$C_x=\RR_{\ge0}$. This shows that $S=M/\approx$ equals $\RR$ with
``doubled'' origin. In particular, it is not Hausdorff.

\def\Problem #1.{\rm\bigbreak\vskip-\parskip\noindent{\bf Problem
#1.}\quad\ignorespaces}

\Problem 6. We have seen that convexity is equivalent to the
combination of \r1, the convexity of $\psi(M)$, \r2, the connectedness
of the fibers of $\psi$, and \r3, the openness of $\psi$ onto its
image.  Properties \r1 and \r2 are well studied in the
literature. Therefore, one may wonder whether \r3 is really an
additional constraint. However, the author was not able to come up
with an example of an Hamiltonian manifold where \r1 and \r2 hold, but
\r3 doesn't. Observe, that Examples~2 and~3 have a non\_convex image
while in the other examples some fibers of $\psi$ are disconnected. In
case, an example with $\r1\wedge\r2\wedge\neg\r3$ exists then the
space $S$ would have have quite peculiar properties: it is a manifold
with non\_empty boundary and possibly corners such that there exists a
continuous {\it bijective\/} map to a {\it convex\/} subset of an
Euclidean space.

\beginsection ConvExamples. Examples of convex Hamiltonian manifolds

Let $V$ be a unitary representation of $K$. Then every smooth
$K$\_stable complex algebraic subvariety of ${\bf P}(V)$ is in a
canonical way a Hamiltonian $K$\_manifold. We call a Hamiltonian
$K$\_manifold {\it projective} if it arises this way but possibly with
the symplectic form and the moment map rescaled by some non\_zero
factor.

\Proposition Proj. Let $M$ be a Hamiltonian $K$\_manifold such that for
every projective Hamiltonian $K$\_manifold $X$ holds that $\psi_{M\times
X}^{-1}(0)$ is connected. Then $\psi_M:M\pfeil\psi(M)$ is an open map
with connected fibers.

\Proof. Let $\Xq$ equal $X$ with the symplectic structure multiplied
by $-1$. Then we have $\mu_{M\times\Xq}(m,x)=\mu_M(m)-\mu_X(x)$. By
choosing for $X$ the coadjoint orbit $Ku$, $u\in\ft^+$ we obtain
that $\psi_M^{-1}(u)=\psi_{M\times\Xq}^{-1}(0)$ is connected. Now let
$X_0$ be projective such that $\psi_{X_0}(X_0)$ is a neighborhood of $0$,
e.g., $X_0={\bf P}(V)$ where $V$ contains stable points for the
$K^\CC$\_action. By rescaling, we can arrange that $\psi(X_0)$ is
arbitrarily small. Let $X:=X_0\times Ku$. Then $B:=\psi_X(X)$ is an
arbitrarily small neighborhood of $u$ in $\ft^+$. Moreover, projection
to the first factor induces a surjective map
$\psi_{M\times\Xq}^{-1}(0)\auf\psi_M^{-1}(B)$. By assumption, the
first, hence the second set, is connected. By
\cite{Conn}, we conclude that $\psi_M$ is open onto its image.\qed

\noindent Now we can give examples of convex Hamiltonian manifolds:

\Theorem Examp. Every connected Hamiltonian $K$\_manifold $M$
satisfying any one of the conditions \r1 through \r4 below is
convex.
\item{\r1} $M$ is compact.
\item{\r2} The moment map $\mu:M\pfeil\fk^*$ is proper.
\item{\r3} The manifold $M$ is a complex algebraic variety, the action
of $K$ is the restriction of an algebraic $K^\CC$\_action, and the symplectic
structure is induced by a $K$\_invariant K{\"a}hler metric.
\item{\r4} The manifold $M$ is a complex Stein space, the action
of $K$ is the restriction of a holomorphic $K^\CC$\_action, and the symplectic
structure is induced by a $K$\_invariant K{\"a}hler metric.
\Par

\Proof. In all cases, it is known that $\psi$ has
connected fibers and convex image (see \cite{HNP} or \cite{Sja} for
\r1, \r2 and \cite{HeHu} for \r3 and \r4). Moreover, the classes
\r1-\r3 are preserved by taking the product with a projective
Hamiltonian $K$\_manifold. Thus \cite{Proj} implies that $\psi$ is
open in these cases.

For \r4 one has to argue a bit differently since the class of Stein
spaces is not stable for products with projective manifolds. But
fortunately, Heinzner\_Huckleberry prove in \cite{HeHu} a stronger
result: given a holomorphic action of $G$ on a complex manifold $M$
then $\psi$ has convex image and connected fibers whenever the
$G$\_action is regular. Here ``regular'' means ``each $G$\_orbit of
$M$ has an open dense neighborhood which is $G$\_equivariantly the
locally biholomorphic image of an irreducible complex $G$\_space $U$
such that $U$ admits a closed, holomorphic, $G$\_equivariant embedding
into an algebraic $G$\_variety.'' It is noted in \cite{HeHu} that
actions on Stein spaces are regular. Clearly, the class of manifolds
with regular action is stable for taking products with algebraic
$G$\_varieties. Thus, $\psi$ is open also in case \r4.\qed

\noindent Note, that \r1 includes all projective spaces while all
linear actions on unitary vector spaces are covered by \r3.

\beginsection LocalConv. Local convexity

\noindent The purpose of this section is to prove the following theorem.

\Theorem SKL. Let $M$ be a Hamiltonian $K$\_manifold. Then every $x\in
M$ has an arbitrarily small $K$\_stable open neighborhood $U$ which is
convex as Hamiltonian manifold and such that $\psi(U)$ is open in
$C_x$.

\Proof. We start with some reductions. Let
$y:=\psi(x)\in\ft^+\subseteq\fk^*$ and $L:=K_y$. Assume $L\ne K$. It
is convenient to identify $\fk^*$ with $\fk$ using a $K$\_invariant
scalar product. Then $y\in\ft\subseteq\fl:=\|Lie|L\subseteq\fk$. Let
$V_0\subseteq\ft^+$ be a small open neighborhood of $y$. Then
$V_L:=L\cdot V_0$ is a small open neighborhood of $Ly$ in
$\fl$. Similarly, $V_K:=K\cdot V_0$ is a small open neighborhood of
$Ky$ in $\fk$. Moreover, $V_L$ is an orthogonal slice to $Ky$. This
implies that $K\times^LV_L\pf\sim V_K$. Now let $M_L:=\mu^{-1}(V_L)$
and $M_K:=\mu^{-1}(V_K)$. Then $M_K$ is an open neighborhood of $Kx$
with $K\times^LM_L\pf\sim M_K$. It can be shown (\cite{GS} Thm.~26.7)
that $M_L$ is a Hamiltonian $L$\_manifold whose moment map $\mu_L$ is
just the restriction of $\mu$ to $M_L$. Similarly, there is $\psi_L$
which is the restriction of $\psi$.  Suppose the theorem is true for
$L$ instead of $K$. Then there is an arbitrarily small open
neighborhood $U_L$ of $Lx$ in $M_L$ which is convex as a Hamiltonian
$L$\_manifold. Then $U:=K\cdot U_L=K\times^LU_L$ is an arbitrarily
small $K$\_invariant open neighborhood of $Kx$. Moreover, for the
preimage of a line segment $\[u,v]\subseteq\ft^*$ holds
$$\eqno{}
\psi^{-1}(\[u,v])\cap U=K\times^L(\psi_L^{-1}(\[u,v])\cap U_L).
$$
This shows that $U$ is a convex Hamiltonian $K$\_manifold.

Thus, it remains to consider the case $L=K$, i.e.,
$y\in(\fk^*)^K$. Since the translate $\mu':=\mu-y$ is again a moment
map we may assume $y=0$. Let $H:=K_x$. Since the orbit $Kx$ is
isotropic we have $\fk x\subseteq(\fk x)^\perp\subseteq T_xM$. Let $S$
be an $H$\_invariant complement of $\fk x$ in $(\fk x)^\perp$. Then
one can show (\cite{Sja} Thm.~6.3, \cite{GS} Thm.~22.1)
\item{$\bullet$}$S$ is a unitary $H$\_module with moment map
$\mu_S:S\pfeil\fh^*:v\mapsto(\xi\mapsto\<v,\xi v\>)$;
\item{$\bullet$}$M_0:=K\times^H(\fh^\perp\oplus S)$ is a Hamiltonian
$K$\_manifold with moment map
$\mu([k,(\alpha,v)])=k(\alpha+\mu_S(v))$.  \item{$\bullet$}Let
$x_0:=[1,(0,0)]\in M_0$. Then $Kx\subseteq M$ and $Kx_0\subseteq M_0$
have $K$\_invariant open neighborhoods which are isomorphic as
Hamiltonian $K$\_manifolds.

\noindent Because of the last point we can (and will) assume $M=M_0$
and $x=x_0$.

Let $D_r$ and $\Dq_r$ be the open and the closed ball of radius $r$ in $S$,
respectively. Consider the natural scalar $U(1)$\_action on $S$ and
let $P=\Dq_r/\sim$ be the space where all $U(1)$\_orbits in
$\Dq_r\setminus D_r$ are collapsed to a point. One can show (Lerman's
symplectic cut, \cite{Le}) that $P$ is a compact Hamiltonian
$H$\_manifold (in fact, it is a complex projective space) whose moment
map $\mu_P$ factors $\mu_S$:
$$\eqno{}
\matrix{\Dq&\auf&P\cr
\dsseq&&\downarrow\Rechts{\mu_P}\cr
S&\pf{\mu_S}&\fh^*\cr}
$$
Now consider
$$\eqno{}
V_r:=K\times^H(\fh^\perp\times D_r),\quad
\Vq_r=K\times^H(\fh^\perp\times \Dq_r),\quad
Q:=K\times^H(\fh^\perp\times P).
$$
Then $V_r$ is an open subset of $M$, its closure is $\Vq_r$, and $Q$
is an image of $\Vq_r$. The point is now that
$\fh^\perp\times\Dq_r\pfeil\fk^*:(\alpha,v)\mapsto\alpha+\mu_S(v)$ is
proper. Therefore, also the moment map of $Q$ is proper. Hence, $Q$ is
convex by \cite{Examp}\r2.

Now we claim that $V_r$ is also convex. First, since $V_r$ is an open
subset of $Q$, we see that the moment map on $V_r$ is open onto its
image. Furthermore, since all fibers of $\Vq_r\pfeil Q$ are connected
(being points or circles) also $\Vq_r\pfeil\ft^+$ has connected
fibers. Now we use that $V_r$ is the union of all $\Vq_s$ with
$s<r$. Thus every fiber of $V_r\pfeil\ft^+$ is an increasing union of
connected sets and therefore connected. The same argument shows that
$\psi(V_r)$ is an increasing union of convex sets, hence itself
convex. This proves the claim.

Let $B_s\subseteq\ft^*$ be the open disk with center $0$ and radius
$s$. From \cite{CP1} we know that $\psi(V_r)$ is a neighborhood of $0$
in $C_x$. Thus, given $r$, we may choose $s$ so small that $B_s\cap
C_x\subseteq\psi(V_r)$. Then $U:=\psi^{-1}(B_s)\cap V_r$ is an
arbitrarily small open neighborhood of $Kx$ which is convex as
Hamiltonian manifold and whose image, $B_s\cap C_x$, is open in $C_x$.
\qed

\Remark. This proof is similar to that of \cite{KL} Prop.~3.7. 
\medskip

\beginsection Property. Some consequences of convexity

If $C$ is a rational cone and $v\in C$ we define the {\it tangent
cone\/} $T_vC$ of $C$ in $v$ as the smallest cone containing $C$ and
having $v$ as a vertex. If $u\in C$ is any vertex then
$$\eqno{2}
T_vC=C+\RR_{\ge0}(u-v).
$$
In fact, to show this we may assume that $u=0$. Then $T_vC$ is a
homogeneous cone containing $C$. The point $v$ being a vertex means
$-v\in T_vC$ which confirms \cite{E2}.

\Corollary CP3. Let $M$ be any Hamiltonian $K$\_manifold and $x\in
M$. Then there is a neighborhood $U$ of $x$ such that for all $y\in U$
we have $C_y=T_{\psi(y)}C_x$.

\Proof. By \cite{SKL} we may assume that $M$ is convex and that
$\psi(M)$ is open in $C_x$. But then for every $u\in\psi(M)$, it holds
that the cone with vertex $u$ spanned by $\psi(M)$ is the same as that
spanned by $C_x$. We conclude with \cite{CP1}.\qed

An application is the following well known statement.

\Corollary fa0. Let $M$ be any connected Hamiltonian $K$\_manifold and
let $\fa\subseteq\ft^*$ be the affine subspace spanned by
$\psi(M)$. Then $\fa$ is also the affine span of every local cone
$C_x$, $x\in M$. Moreover, the set of interior points of $\psi(M)$
(relative to $\fa$) is dense in $\psi(M)$.

\Proof. Let $\fa_x$ be the affine span of $C_x$. It follows from
\cite{CP3} that is locally constant in $x$. Since $M$ is connected,
$\fa_x$ is independent of $x$. From $\psi(x)\in C_x\subseteq\fa_x$
follows $\psi(M)\subseteq\fa_x$. Moreover, for every $x\in M$ there is
a subset $U_x\subseteq\psi(M)\cap C_x$ which is open in $\fa_x$ and
which contains $\psi(x)$ in its closure (\cite{SjaThm}). Thus $\fa_x$
is spanned by $\psi(M)$ which implies $\fa_x=\fa$. Moreover, $\psi(x)$
is in the closure of the relative interior points of $\psi(M)$.\qed

The last class of examples for convex Hamiltonian manifolds was
suggested to me by Tolman:

\Theorem. Let $X$ be a connected $K$\_manifold and let $M:=T^*(X)$ with
its natural Hamiltonian structure. Then $M$ is convex.

\Proof. Let $\iota:X\into M$ be the zero\_section. The
map $\psi:M\pfeil\ft^*$ is homogeneous with respect to the natural
scalar $\RR_{>0}$\_action on the fibers. Therefore,
$\iota(X)\subseteq\psi^{-1}(0)$. Since $X$ is connected, \cite{SjaThm}
implies that $C=C_{\iota(x)}$ is independent of $x\in X$. Moreover,
there is an open neighborhood $U$ of $\iota(X)$ in $M$ such that
$\psi(U)$ is a neighborhood of $0$ in $C$. From $M=\cup_{t>0}tU$ we
obtain $\psi(M)=C$. In particular, we see that $\psi(M)$ is convex.

Let $\pi:M\pfeil X$ be the projection. For $m\in M$ let
$m_0:=\iota\pi(m)\in M$. Choose a convex open neighborhood $V$ of
$Km_0$ such that $\psi(V)$ is open in $C$. Using the
$\RR_{>0}$\_action we may assume that $m\in V$. Now let $V_0$ be a
small neighborhood of $m$ which is contained in $V$. Using the
convexity of $V$, we conclude that $\psi(V_0)$ is open in $\psi(V)$,
hence open in $C$. Thus, $\psi:M\pfeil C$ is open.

Finally, suppose that the fiber of $\psi$ over $u\in C$ is
disconnected. Thus, $ \psi^{-1}(u)=F_1\cup F_2$ where $F_1$ and $F_2$
are disjoint, non\_empty, and closed. Let $W\subseteq X$ be open,
non\_empty and $K$\_stable. Then, with $X$ replaced by $W$, we obtain
as above $\psi(\pi^{-1}(W))=\psi(T^*(W))=C$. This means in particular
that $\pi^{-1}(W)\cap\psi^{-1}(u)\ne\emptyset$ for every $W$.  Hence,
$\pi(\psi^{-1}(u))$ is dense in $X$. Since $X$ is connected there is
$x\in\overline{\pi(F_1)}\cap\overline{\pi(F_2)}$. Let $V$ be a convex
neighborhood of $K\iota(x)$ in $M$. By homogeneity, we may assume that
$V$ meets both $F_1$ and $F_2$. But this contradicts the connectedness
of fibers of $\psi|_V$. Thus, $\psi$ has connected fibers as well.\qed

\Remark. The fact that $\psi(M)=C_x$ for any $x\in X\subseteq T^*(X)$
is due to Sjamaar, \cite{Sja}~Thm.~7.6.

\let\[=\oldbracket

\lastpage
\bye